\documentclass[12pt,a4paper,twoside]{article}
\usepackage{amsfonts}
\date{}
\newcommand{\A}[1]{\vspace{3mm}}
\newcommand{\D}{{\rm d}}

\newcommand{\R}{{\mathbb R}}

\def\S{{\mathbb S}}
\def\H{{\mathbb H}}

\begin{document}

\vspace*{30mm}

\begin{center}
{\Large \bf Covering Large Balls with Convex Sets\\[3mm] in Spherical Space}\\[7mm]
{\bf K\'{a}roly Bezdek \quad Rolf Schneider}

\vspace{3mm}

{\em
Department of Mathematics and Statistics, University of Calgary\\
2500 University Drive N.W., AB, Canada, T2N 1N4\\
e-mail: bezdek@math.ucalgary.ca

\vspace{3mm}

Mathematisches Institut, Albert-Ludwigs-Universit\"at\\
Eckerstr.~1, D-79104 Freiburg i.~Br., Germany\\
e-mail: rolf.schneider@math.uni-freiburg.de}

\end{center}

\vspace{10mm}

\begin{quote}
{\bf Abstract.} If the $n$-dimensional unit sphere is covered by finitely many spherically convex bodies, then the sum of the inradii of these bodies is at least $\pi$. This bound is sharp, and the equality case is characterized.\\[2mm]
\end{quote}

\renewcommand{\thefootnote}{{}}

\footnote{Partially supported by a Natural Sciences and Engineering Research Council of Canada Discovery Grant. R. Schneider thanks the PIMS `Collaborative Research Group in Geometric and Harmonic Analysis' for an invitation to the universities of Edmonton and Calgary, where this note grew out of discussions. }

{\bf 1. Introduction}

Let $\S^n$ be the $n$-dimensional unit sphere in $(n+1)$-dimensional Euclidean space $\R^{n+1}$ ($n\ge 2$). A {\em spherically convex body} is a closed, spherically convex subset $K$ of $\S^n$ with interior points and lying in some closed hemisphere, thus, the intersection of $\S^n$ with an $(n+1)$-dimensional closed convex cone of $\R^{n+1}$ different from $\R^{n+1}$. The {\em inradius} $r(K)$ of $K$ is the spherical radius of the largest spherical ball contained in $K$.

The purpose of this note is the proof of the following theorem.

\vspace{3mm}

{\bf Theorem 1.} {\em If the spherically convex bodies $K_1,\dots,K_m$ cover the spherical ball $B$ of radius $r(B) \ge\pi/2$ in $\S^n$, then
$$ \sum_{i=1}^m r(K_i)\ge r(B).$$
For $r(B)=\pi/2$ the stronger inequality $\sum_{i=1}^m r(K_i\cap B)\ge r(B)$ holds. Moreover, equality for $r(B)=\pi$ or $r(B)=\pi/2$ holds if and only if $K_1,\dots,K_m$ are lunes with common ridge which have pairwise no common interior points.}

\vspace{3mm}

Recall that a {\em lune} in $\S^n$ is the $n$-dimensional intersection of $\S^n$ with two closed half\-spaces of $\R^{n+1}$ with the origin $0$ in their boundaries. The intersection of the boundaries (or any $(n-1)$-dimensional subspace in that intersection, if the two subspaces are identical) is called {\em ridge} of the lune. Evidently, the inradius of a lune is half the interior angle between the two defining hyperplanes.

The original motivation for Theorem 1 came from Tarski's plank problem. It states that if a convex body in Euclidean space is covered by finitely many slabs, then the sum of their widths is at least the (minimal) width of the body. The problem was solved by Bang \cite{Ba50, Ba51}. 
Several related questions are discussed in \cite{Gar88} and \cite[Section 3.4]{BMP05}. The symmetric case of a more general conjecture of Bang was proved by Ball \cite{B91}, who considered coverings of balls by planks in finite-dimensional Banach spaces (the width of a plank being defined in terms of the norm). Recently the plank theorem was further strengthened by Kadets \cite{Ka05} for Hilbert spaces, as follows. Let $C$ be a convex body, i.e., a closed convex subset with non-empty interior, in the real Hilbert space $\H$ (finite- or infinite-dimensional). Let $r(C)$ denote the supremum of the radii of the balls contained in $C$. Planks and their widths are defined with the help of the inner product of $\H$ in the usual way. Thus, if $C$ is a convex body in $\H$ and $P$ is a plank of $\H$, then the width of $P$ is always at least as large as $2r(C\cap P)$. The result of \cite{Ka05} says that if a convex body $B$ (it suffices to consider balls) of $\H$ is covered by the convex bodies $C_1,\dots , C_n$ in $\H$, then $\sum_{i=1}^n r(C_i)\ge \sum_{i=1}^n r(C_i\cap B)\ge r(B)$. We note that an independent proof of the $2$-dimensional Euclidean case of this result can be found in \cite{Bez07}. The proofs of \cite{Ka05} and \cite{Bez07} do not generalize neither to Banach spaces nor to spherical space. While an extension to Banach spaces seems rather difficult, we noticed that in spherical space for coverings of large balls an easier answer can be given.

Theorem 1 is a consequence of the following result. Here $\sigma$ denotes spherical Lebesgue measure on $\S^n$, and $\sigma_n:= \sigma(\S^n )$.

\vspace{3mm}

{\bf Theorem 2.} {\em If $K$ is a spherically convex body, then
$$ \sigma(K)\le \frac{\sigma_n}{\pi} r(K).$$
Equality holds if and only if $K$ is a lune.}

\vspace{3mm}

This implies Theorem 1 as follows. If $B=\S^n$, i.e., the spherically convex bodies $K_1,\dots,K_m$ cover $\S^n$, then
$$ \sigma_n \le \sum_{i=1}^m \sigma(K_i) \le \frac{\sigma_n}{\pi} \sum_{i=1}^m r(K_i),$$
and the stated inequality follows. In general, when $B$ is different from $\S^n$,
let $B'\subset \S^n$ be the spherical ball of radius $\pi-r(B)$ centered at the point antipodal to the center of $B$.
As the spherically convex bodies $B', K_1,\dots,K_m$ cover $\S^n$, the inequality just proved shows that
$$ \pi-r(B)+\sum_{i=1}^m r(K_i)\ge \pi,$$
and the stated inequality follows. If $r(B)=\pi/2$, then $K_1\cap B,\dots,K_m\cap B$ are spherically convex bodies and as $B', K_1\cap B,\dots,K_m\cap B$ cover $\S^n$ , the stronger inequality follows. The assertion about the equality sign for the case when $r(B)= \pi$ or $\pi/2$ follows easily.

\vspace{5mm}

{\bf 2. Proof of Theorem 2}

We denote the standard scalar product of $\R^{n+1}$ by $\langle\cdot,\cdot\rangle$, and for $u\in\S^n$ we write
$$u^\perp:=\{x\in\R^{n+1}: \langle u,x \rangle =0\}$$
for the orthogonal complement of ${\rm lin}\{u\}$. For a spherically convex body $K$, the polar body is defined by
$$ K^*:=\{u\in \S^n: \langle u,v \rangle \le 0 \mbox{ for all }v\in K\}.$$
It is also spherically convex, but need not have interior points. The number
$$ U(K) :=\frac{1}{2}\sigma(\{u\in\S^n:u^\perp\cap K\not=\emptyset\})$$
can be considered as the {\em spherical mean width} of $K$. Obviously, a vector $u\in\S^n$ satisfies $u\in K^*\cup(-K^*)$ if and only if $u^\perp$ does not meet the interior of $K$, hence
\begin{equation}\label{2.1}
\sigma_n-2\sigma(K^*) = 2U(K).
\end{equation}
It is a basic idea of the following to treat $U(K^*)$ instead of $\sigma(K)$.

Let $K\subset \S^n$ be a spherically convex body, and let $B$ be the smallest spherical ball containing it. We assume that $B$ is not a closed hemisphere (this will be satisfied later). Let $e\in\S^n$ be the center of $B$, let $\S^n_e:=\{u\in \S^n: \langle e,u\rangle>0\}$, and let $T_e$ be the tangent hyperplane to $\S^n$ at $e$. We write ${\mathbb E}^n:=e^\perp$. With the induced scalar product, this is an $n$-dimensional Euclidean space. Further, we write $\S^{n-1}:= \S^n\cap{\mathbb E}^n$ for the unit sphere of this space. The mapping $\Pi:\S^n_e\to {\mathbb E}^n$ is defined as the radial projection from $\S^n_e$ to $T_e$, followed by the orthogonal projection to ${\mathbb E}^n$.

Every $(n-1)$-dimensional great subsphere of $\S^n$ which does not contain $e$ can uniquely be written in the form $S_u = u^\perp \cap \S^n$ with $u\in \S_e^n$. For $u\in\S_e^n\setminus\{e\}$ the subspace $u^\perp$ intersects the hyperplane $T_e$ in an $(n-1)$-dimensional affine subspace, and under orthogonal projection to ${\mathbb E}^n$ this is projected into a hyperplane of ${\mathbb E}^n$. Writing
$$ u= \tau e-\sqrt{1-\tau^2}\,u_0, \qquad 0\le \tau < 1, \quad u_0\in\S^{n-1},$$
we find that this hyperplane is given by
$$ \Pi(S_u \cap \S^n_e) = H(u_0,t) :=\{x\in {\mathbb E}^n: \langle u_0,x \rangle =t\}$$
with
$$ t=t(\tau)= \frac{\tau}{\sqrt{1-\tau^2}}.$$

The set $\overline{K}:=\Pi(K)$ is a compact convex set in ${\mathbb E}^n$, and $\overline{B}:=\Pi(B)$ is the smallest Euclidean ball containing $\overline{K}$. The ball $\overline{B}$ has center $0$. Since $u^\perp\cap K\not=\emptyset$ is equivalent to $H(u_0,t(\tau))\cap \overline{K}\not=\emptyset$, we obtain (noting that $\sigma(\{u\in \S^n: e\in u^\perp\})=0$)
\begin{eqnarray*}
U(K) &=& \int_{\S^n_e} {\bf 1}\{u^\perp \cap K\not=\emptyset\}\,\sigma(\D u)\\
&=& \int_{\S^{n-1}} \int_0^1 {\bf 1}\{H(u_0,t(\tau)) \cap \overline{K}\not=\emptyset\}(1-\tau^2)^{\frac{n-2}{2}}\,\D\tau\,\mu(\D u_0),
\end{eqnarray*}
where $\mu$ denotes spherical Lebesgue measure on $\S^{n-1}$. Substituting $\tau = t/\sqrt{1+t^2}$, we obtain
$$ U(K) = \int_{\S^{n-1}} \int_0^\infty {\bf 1}\{H(u_0,t) \cap \overline{K}\not=\emptyset\}(1+t^2)^{-\frac{n+1}{2}}\,\D t\,\mu(\D u_0).$$

We can now state a more general assertion, from which Theorem 2 will follow. For this, let $f:[0,\infty)\to (0,\infty)$ be a positive continuous function. On the space ${\mathcal H}$ of hyperplanes of ${\mathbb E}^n$ (with its usual topology) we define a Borel measure $\nu_f$ by
$$ \nu_f(A):= \int_{\S^{n-1}} \int_0^\infty {\bf 1}\{H(u,t) \in A\}f(t)\,\D t\,\mu(\D u)$$
for Borel sets $A\subset {\mathcal H}$. Then
$$ U_f(\overline{K}) := \int_{\mathcal H} {\bf 1}\{H\cap \overline{K}\not=\emptyset\}\,\nu_f(\D H) $$
is the total $\nu_f$ measure of the set of hyperplanes meeting the convex set $\overline K\subset {\mathbb E}^n$, and
$$ U(K) = U_f(\overline{K}) \qquad \mbox{for } f(t)= (1+t^2)^{-\frac{n+1}{2}}. $$

\vspace{3mm}

\noindent{\bf Proposition.} {\em Let $B\subset {\mathbb E}^n$ be a ball with center $0$, and let ${\mathcal K}(B)$ be the set of nonempty, compact convex sets $K\subset {\mathbb E}^n$ for which $B$ is the smallest ball containing $K$. Then, for $K\in{\mathcal K}(B)$, the value $U_f(K)$ is minimal if and only if $K$ is a segment.}

\vspace{3mm}

For a constant function $f$, in which case $U_f$ is a constant multiple of the mean width, this result is due to Linhart \cite{Lin77}. We extend Linhart's proof, but also modify it, for a reason explained later.

Let $K\in{\mathcal K}(B)$. Writing $h(K,u):=\max\{\langle x,u\rangle: x\in K\}$ for the support function of $K$ at $u\in\S^{n-1}$, we have (noting that $0\in K$)
\begin{eqnarray*}
U_f(K) &=& \int_{\S^{n-1}} \int_0^\infty {\bf 1}\{H(u,t) \cap K\not=\emptyset\}f(t)\,\D t\,\mu(\D u)\\
&=& \int_{\S^{n-1}} \int_0^{h(K,u)} f(t)\,\D t\,\mu(\D u) = \int_{\S^{n-1}} F(h(K,u))\,\mu(\D u)
\end{eqnarray*}
with $F(s):= \int_0^s f(t)\,\D t$.

By a standard compactness and continuity argument, $U_f$ attains a minimum on ${\mathcal K}(B)$, say at $K$. Since $K\in{\mathcal K}(B)$, there exists a $k$-simplex $T$ with vertices $v_1,\dots,v_{k+1} \in {\rm bd}\,B$, for some $k\in\{1,\dots,n\}$, such that $T\subset K$ and $T\in{\mathcal K}(B)$. Then $K=T$, since otherwise there is an open set of hyperplanes of positive $\nu_f$ measure hitting $K$ but not $T$, which would imply $U_f(T)< U_f(K)$. For $j=1,\dots,k+1$ let $N(T,v_j)$ denote the normal cone of $T$ at its vertex $v_j$, so that $S_j:= N(T,v_j)\cap \S^{n-1}$ is the spherical image of $v_j$. For $u\in S_j$ we have $h(T,u)= \langle v_j,u\rangle =R\cos\varphi$, where $R$ is the radius of $B$ and $\varphi= \varphi(u)$ denotes the angle between $v_j$ and $u$. Let $D_j:= \{u\in\S^{n-1}: \langle u,v_j\rangle \ge 0\}$; then $S_j\subset D_j$. We write $g(\varphi) := F(R\cos\varphi)$ and state a generalization of Linhart's crucial inequality, namely
\begin{equation}\label{7.1}
\frac{1}{\mu(S_j)} \int_{S_j} g(\varphi(u))\,\mu(\D u)
\ge \frac{1}{\mu(D_j)} \int_{D_j} g(\varphi(u))\,\mu(\D u) =: C(R,f),
\end{equation}
with equality if and only if $S_j=D_j$. On the right side, $C(R,f)$ denotes a constant that depends on $R$ and on the function $f$, but is independent of $j$. If (\ref{7.1}) has been proved, then it follows that
\begin{eqnarray*}
U_f(K)=U_f(T) &=& \int_{\S^{n-1}} F(h(T,u))\,\mu(\D u) = \sum_{j=1}^{k+1} \int_{S_j} F(h(T,u))\,\mu(\D u) \\
&\ge& \sum_{j=1}^{k+1} \mu(S_j)C(R,f) = \mu(\S^{n-1})C(R,f),
\end{eqnarray*}
with equality if and only if $T$ is a segment.

To prove (\ref{7.1}), we modify the proof given by Linhart, replacing the function $R\cos\varphi$ by the function $g(\varphi)$.
Moreover, instead of approximation by step functions, we use an integration argument. The reason for this is that the use of approximation blurs the equality cases, so that we would not be able to conclude that the minimum is attained only by segments.

We fix $j\in\{1,\dots,k+1\}$ and write $e_j:= v_j/R$ and, for $u\in S_j$,
$$ u=e_j\cos\varphi+w\sin\varphi,\qquad 0\le\varphi\le\pi/2,\quad w\in W:= \S^{n-1}\cap e_j^\perp, $$
so that $\varphi=\varphi(u)=\langle e_j,u\rangle$. For $w\in W$, there is a unique value $b=b(w)$ such that $e_j\cos \varphi+w\sin \varphi \in S_j$ precisely for $0\le \varphi\le b$. We write $J(\varphi):=\sin^{n-2}\varphi$ and denote by $\D w$ the $(n-2)$-dimensional spherical surface area element of the sphere $W$ at $w$; then
$$ \int_W \int_0^{\pi/2}J(\varphi)\,\D\varphi\,\D w = \mu(D_j). $$
Using Fubini's theorem, we can write
\begin{eqnarray*}
\int_{S_j} g(\varphi(u))\,\mu(\D u) &=& \int_W \int_0^{b(w)} g(\varphi)J(\varphi)\,\D\varphi\,\D w\\
&=& \int_0^\infty\int_W \int_0^{b(w)} {\bf 1}\{t\le g(\varphi)\}J(\varphi)\,\D\varphi\,\D w\,\D t.
\end{eqnarray*}
Let $t\ge 0$ be fixed. With
$$ M(t) := \int_W \int_0^{\pi/2} {\bf 1}\{t\le g(\varphi)\}J(\varphi)\,\D\varphi\,\D w$$
we have
$$ \frac{M(t)}{\mu(D_j)} = \frac{\int_0^{\pi/2} {\bf 1}\{t\le g(\varphi)\}J(\varphi)\,\D \varphi}{\int_0^{\pi/2}J(\varphi)\,\D \varphi} \le 1.$$

Let $A_1:= \{w\in W:g(b(w))\ge t\}$ and $A_2:= \{w\in W:g(b(w))< t\}$. Since the function $F$ is strictly increasing and $\cos $ is strictly decreasing on $[0,\pi/2]$, the function $g$ is strictly decreasing.
For $w\in A_1$ and $\varphi\le b(w)$ we have $g(\varphi)\ge g(b(w))\ge t$, hence ${\bf 1}\{t\le g(\varphi)\}=1$.
For $w\in A_2$ and $\varphi>b(w)$ we have $g(\varphi)<g(b(w))< t$, hence ${\bf 1}\{t\le g(\varphi)\}=0$. This gives
\begin{eqnarray*}
& & \int_W \int_0^{b(w)} {\bf 1}\{t\le g(\varphi)\}J(\varphi)\,\D \varphi\,\D w\\
& & = \int_{A_1} \int_0^{b(w)} J(\varphi)\,\D \varphi\,\D w + \int_{A_2}\int_0^{\pi/2} {\bf 1}\{t\le g(\varphi)\}J(\varphi)\,\D \varphi\,\D w\\
& & = \int_{A_1} \int_0^{b(w)} J(\varphi)\,\D\varphi\,\D w + \frac{M(t)}{\mu(D_j)}\int_{A_2}\int_0^{\pi/2} J(\varphi)\,\D\varphi\,\D w \\
& & \ge \frac{M(t)}{\mu(D_j)} \int_W \int_0^{b(w)} J(\varphi)\,\D \varphi\,\D w = \frac{M(t)}{\mu(D_j)}\mu(S_j).
\end{eqnarray*}
Now an integration over $t$ yields the assertion (\ref{7.1}). Equality holds if and only if $b(w)=\pi/2$ for all $w\in W$, that is, if $S_j=D_j$. If this holds for $j=1,\dots,k+1$, then $k=1$, hence $K$ is a segment. This completes the proof of the Proposition.

To deduce the assertion of Theorem 2, we assume that $K\subset\S^n$ is a spherically convex body with prescribed inradius $r$. The smallest spherical ball containing the polar body $K^*$ has radius $\pi/2-r<\pi/2$. From the Proposition we deduce that $U(K^*)$ is minimal if and only if $\Pi(K^*)$ is a segment, hence if and only if $K$ is a lune. Now it follows from (\ref{2.1}) (with $K$ and $K^*$ interchanged) that the spherical volume $\sigma(K)$ is maximal if and only if $K$ is a lune. The spherical volume of a lune of inradius $r$ is given by $(\sigma_n/\pi)r$.

\end{document}